\documentclass[11pt,reqno]{amsart}
\usepackage{amsmath,amsfonts, amssymb, amsthm}
  \usepackage{paralist}
  \usepackage{graphics} 
  \usepackage{epsfig} 
\usepackage{graphicx}  \usepackage{epstopdf}
 \usepackage[colorlinks=true]{hyperref}
\hypersetup{urlcolor=blue, citecolor=red}

\newtheorem{theorem}{Theorem}[section]

\newtheorem{lemma}[theorem]{Lemma}

\theoremstyle{definition}
\newtheorem{definition}{Definition}[section]
\newtheorem{remark}{Remark}[section]

\def\pmod #1{\ ({\rm{mod}}\ #1)}
\def\Z{\Bbb Z}
\def\N{\Bbb N}

\def\Q{\Bbb Q}

\def\R{\Bbb R}

\def\l{\left}
\def\r{\right}
\def\bg{\bigg}
\def\({\bg(}
\def\){\bg)}
\def\t{\text}
\def\f{\frac}
\def\mo{{\rm{mod}\ }}
\def\pmod#1{\ (\mo\ #1)}

\def\ls{\leq}
\def\gs{\geq}
\def\se {\subseteq}
\def\sm{\setminus}

\def\al{\alpha}

\def\ve{\varepsilon}

\def\eq{\equiv}

\def\Proof{\noindent{\it Proof}}
\def\Ack{\medskip\noindent {\bf Acknowledgment}}

\begin{document}
\hbox{Bull. Pol. Acad. Sci. Math. 70 (2022), no.\,2, 93--106.}
\medskip

\title[$\mathbb Q\setminus\mathbb Z$ is diophantine over $\mathbb Q$ with $32$ unknowns]
      {$\Q\sm\Z$ is diophantine over $\Q$ with $32$ unknowns}
\author[Geng-Rui Zhang and Zhi-Wei Sun]{Geng-Rui Zhang and Zhi-Wei Sun$^\star$}


\address{(Geng-Rui Zhang) School of Mathematical Sciences, Peking
University, Beijing 100871, People's Republic of China}
\email{grzhang@stu.pku.edu.cn}

\address{(Zhi-Wei Sun, corresponding author) Department of Mathematics, Nanjing
University, Nanjing 210093, People's Republic of China}
\email{zwsun@nju.edu.cn}

\keywords{Undecidability, definability, diophantine sets, Hilbert's tenth problem over $\Q$, mixed quantifiers.
\newline \indent
$^\star$ Corresponding author, supported by the National Natural Science Foundation of China (grant no. 11971222).}

\begin{abstract} In 2016 J. Koenigsmann refined a celebrated theorem of J. Robinson by proving that $\mathbb Q\setminus\mathbb Z$ is diophantine over $\Q$, i.e., there is a polynomial $P(t,x_1,\ldots,x_{n})\in\Z[t,x_1,\ldots,x_{n}]$
such that for any rational number $t$ we have
$$t\not\in\Z\iff \exists x_1\cdots\exists x_{n}[P(t,x_1,\ldots,x_{n})=0]$$
where variables range over $\mathbb Q$, equivalently
$$t\in\Z\iff \forall x_1\cdots\forall x_{n}[P(t,x_1,\ldots,x_{n})\not=0].$$
In this paper we prove that we may take $n=32$.
Combining this with a result of Z.-W. Sun, we show that there is no algorithm to decide for any $f(x_1,\ldots,x_{41})\in\Z[x_1,\ldots,x_{41}]$ whether
$$\forall x_1\cdots\forall x_9\exists y_1\cdots\exists y_{32}[f(x_1,\ldots,x_9,y_1,\ldots,y_{32})=0],$$
where variables range over $\mathbb Q$.
\end{abstract}
\maketitle

\section{Introduction}

Hilbert's Tenth Problem (HTP) asks for an algorithm to determine for any given polynomial $P(x_1,\ldots,x_n)\in\Z[x_1,\ldots,x_n]$
whether the diophantine equation $P(x_1,\ldots,x_n)=0$ has solutions $x_1,\ldots,x_n\in\Z$.
This was solved negatively by Yu. Matiyasevich \cite{M70} in 1970, on the basis of the important work
of M. Davis, H. Putnam and J. Robinson \cite{DPR}; see also Davis \cite{D73} for a nice introduction.
Z.-W. Sun \cite{S21} proved his 11 unknowns theorem which states that there is no algorithm to determine for any $P(x_1,\ldots,x_{11})\in\Z[x_1,\ldots,x_{11}]$ whether the equation $P(x_1,\ldots,x_{11})=0$
has solutions over $\Z$.

It remains open whether HTP over $\Q$ is undecidable.
However, Robinson \cite{R49} used the theory of quadratic forms to prove that one can characterize
$\Z$ by using the language of $\Q$ in the following way: For any $t\in\Q$ we have
$$t\in\Z\iff\forall x_1\forall x_2\exists y_1\cdots\exists y_7\forall z_1\cdots\forall z_6[f(t,x_1,x_2,y_1,\ldots,y_7,z_1,\ldots,z_6)=0],$$
where $f$ is a polynomial with integer coefficients. (Throughout this paper, variables always range over $\Q$.) In 2009 B. Poonen \cite{Po} improved this by finding a polynomial $F(t,x_1,x_2,y_1,\ldots,y_7)$ with integer coefficients such that for any $t\in\Q$ we have
$$t\in \Z\iff \forall x_1\forall x_2\exists y_1\cdots\exists y_7[F(t,x_1,x_2,y_1,\ldots,y_7)=0].$$
In 2016 J. Koenigsmann \cite{K} improved Poonen's result by proving that the set $\mathbb Q\setminus\mathbb Z$ is diophantine over $\Q$, i.e., there is a polynomial $P(t,x_1,\ldots,x_{n})\in\Q[t,x_1,\ldots,x_{n}]$
such that for any $t\in\Q$ we have
$$t\not\in\Z\iff \exists x_1\cdots\exists x_{n}[P(t,x_1,\ldots,x_{n})=0],$$
i.e.,
$$t\in\Z\iff \forall x_1\cdots\forall x_{n}[P(t,x_1,\ldots,x_{n})\not=0].$$
The number $n$ of unknowns in Koenigsmann's diophantine representation of $\Q\sm \Z$ over $\Q$
is over 400 but below 500. In 2018 N. Daans \cite{D18} significantly simplified Koenigsmann's
approach and proved that $\Q\sm\Z$ has a diophantine representation over $\Q$
with $50$ unknowns. The number $50$ could be reduced to $38$
by applying a recent result \cite[Theorem 1.4]{DDF} obtained by model theory.

In this paper we establish the following new result.

\begin{theorem} \label{Th1.1} $\Q\sm \Z$ has a diophantine representation over $\Q$ with $32$ unknowns, i.e., there is a polynomial $P(t,x_1,\ldots,x_{32})\in\Z[t,x_1,\ldots,x_{32}]$ such that for any $t\in\Q$ we have
\begin{equation}\label{32} t\not\in\Z\iff \exists x_1\cdots\exists x_{32}[P(t,x_1,\ldots,x_{32})=0].
\end{equation}
Furthermore, the polynomial $P$ can be constructed explicitly with $\deg P<2.1\times10^{11}$.
\end{theorem}

To obtain this theorem, we start from Daans' work \cite{D18}, and mainly use a new relation-combining theorem on diophantine representations over $\Q$ (which is an analogue of Matiyasevich and Robinson's relation-combining theorem \cite[Theorem 1]{MR})
as an auxiliary tool.  Now we state our relation-combining theorem for diophantine representations over $\Q$.

\begin{theorem}\label{Th1.2} Let ${\mathcal J}_k(x_1,\ldots,x_k,x)$ denote the expression
\begin{align*}\prod_{s=1}^kx_s^{(k-1)2^{k+1}}
\times\prod_{\ve_1,\ldots,\ve_k\in\{\pm1\}}
\l(x+\sum_{s=1}^k\ve_s\sqrt{x_s}\,W(x_1,\ldots,x_k)^{s-1}\r),
\end{align*}
where
$$W(x_1,\ldots,x_k)=\l(k+\sum_{s=1}^kx_s^2\r)\l(1+\sum_{s=1}^kx_s^{-2}\r).$$
Then ${\mathcal J}_k(x_1,\ldots,x_k,x)$ is a polynomial with integer coefficients. Moreover,
for any $A_1,\ldots,A_k\in\Q^*=\Q\sm\{0\}$, we have
\begin{equation}A_1,\ldots,A_k\in\square \iff \exists x [{\mathcal J}_k(A_1,\ldots,A_k,x)=0],
\end{equation}
where $\square=\{r^2:\ r\in\Q\}$.
\end{theorem}
\begin{remark} In view of its proof, Theorem \ref{Th1.2} can be generalized by replacing $\Q$ with any subfield of the real field $\R$ or any ordered field.
\end{remark}

When $\rho_s\in\{\forall,\exists\}$ for all $s=1,\ldots,k$, we say that $\rho_1\cdots\rho_k$
over $\Q$ is undecidable if there is no algorithm to decide for any polynomial $P(x_1,\ldots,x_k)$ over $\Q$
whether
$$\rho_1 x_1\cdots \rho_k x_k[P(x_1,\ldots,x_k)=0]$$
or not. For convenience we adopt certain abbreviation, for example, $\forall^2\exists^3$
denotes $\forall\forall \exists\exists\exists$.

Combining Theorem \ref{Th1.1} and its proof with a result of Sun \cite[Theorem 1.1]{S21}, we obtain the following theorem.

\begin{theorem}\label{Th1.3}  $\forall^{9}\exists^{32}$ over $\Q$ is undecidable,
 i.e., there is no algorithm
to determine for any $P(x_1,\ldots,x_{41})\in\Z[x_1,\ldots,x_{41}]$ whether
$$\forall x_1\cdots\forall x_{9}\exists y_1\cdots\exists y_{32}[P(x_1,\ldots,x_{9},y_1,\ldots,y_{32})=0].$$
Also,
$\exists^9\forall^{32}\exists$ over $\Q$ and $\exists^{10}\forall^{31}\exists$
over $\Q$ are undecidable.
\end{theorem}

We remark that Sun \cite{Sun} obtained some undecidability results on mixed quantifier prefixes over diophantine equations with integer variables; for example, he proved that $\forall^2\exists^4$ over $\Z$ is undecidable.

In the next section we will prove Theorem \ref{Th1.2}.
Sections 3 and 4 are devoted to our proofs of Theorems \ref{Th1.1} and \ref{Th1.3} respectively.

\section{Proof of Theorem \ref{Th1.2}}
 \setcounter{equation}{0}
 \setcounter{conjecture}{0}
 \setcounter{theorem}{0}

  \medskip
  \noindent{\it Proof of Theorem \ref{Th1.2}}. Clearly,
  $$I_k(x_1,\ldots,x_k,x,y)=\prod_{\ve_1,\ldots,\ve_k\in\{\pm1\}}(x+\ve_1x_1+\ve_2x_2y+\cdots+\ve_kx_ky^{k-1}).$$
  is a polynomial with integer coefficients. As
  $$I_k(x_1,\ldots,x_k,x,y)=\prod_{\ve_i\in\{\pm1\}\ \t{for}\ i\not=t}\(\(x+\sum_{s=1\atop s\not=t}^k\ve_sx_sy^{s-1}\)^2-x_t^2y^{2(t-1)} \)$$
  for all $t=1,\ldots,k$, we see that
  $$I_k(x_1,\ldots,x_k,x,y)=I_k^*(x_1^2,\ldots,x_k^2,x,y)$$
  for some polynomial $I_k^*$ with integer coefficients. Note that
  \begin{align*}{\mathcal J}_k(x_1,\ldots,x_k,x)=&\prod_{s=1}^kx_s^{(k-1)2^{k+1}}
  \\&\times I_k^*\(x_1,\ldots,x_k,x,\(k+\sum_{j=1}^kx_j^2\)\(1+\sum_{j=1}^kx_j^{-2}\)\)
  \end{align*}
  is a polynomial with integer coefficients.

  Now let $A_1,\ldots,A_k\in\Q^*$. We claim that for any rational number
  \begin{equation}
  \label{W-bound} W_k\gs\f{1+\sum_{s=1}^k|\sqrt{A_s}|}{\min\{|\sqrt{A_1}|,\ldots,|\sqrt{A_k}|\}},
  \end{equation}
  we have
  $$A_1,\ldots,A_k\in\square\iff \exists x[I_k^*(A_1,\ldots,A_k,x,W_k)=0].$$

   The ``$\Rightarrow$" direction is easy. If $A_1=a_1^2,\ldots,A_k=a_k^2$ for some $a_1,\ldots,a_k\in\Q$, then, for $x=a_1+a_2W_k+\cdots+a_kW_k^{k-1}\in \Q$ we  have
  $I_k^*(A_1,\ldots,A_k,x,W_k)=0$.

  We use induction on $k$ to prove the ``$\Leftarrow$" direction of the claim.
  In the case $k=1$, if $I_1^*(A_1,x,W_1)=x^2-A_1$ is zero for some $x\in\Q$ then we obviously have $A_1\in\square$.

  Now let $k>1$ and assume that the ``$\Leftarrow$" direction of the claim
  holds for all smaller values of $k$. Let $W_k$ be any rational number satisfying the inequality \eqref{W-bound}.
  Suppose that $I_k^*(A_1,\ldots,A_k,x,W_k)=0$ for some $x\in\Q$. Then there are $\ve_1,\ldots,\ve_k\in\{\pm1\}$ such that
  $$x+\sum_{s=1}^k\ve_s\sqrt{A_s}W_k^{s-1}=0.$$
  If $A_k=a_k^2$ for some $a_k\in\Q$, then, for $x'=x+\ve_k|a_k|W_k^{k-1}$ we have
  $$x'+\ve_1\sqrt{A_1}+\ve_2\sqrt{A_2}W_k+\cdots+\ve_{k-1}\sqrt{A_{k-1}}W_k^{k-2}=0$$
  and hence $I_{k-1}^*(A_1,\ldots,A_{k-1},x',W_k)=0$. Note that
  $$|\sqrt{A_t}|W_k\gs1+\sum_{s=1}^k|\sqrt{A_s}|\gs1+\sum_{s=1}^{k-1}|\sqrt{A_s}|$$
   for each $t=1,\ldots,k-1$. So, in the case $A_k\in\square$, we get $A_1,\ldots,A_{k-1}\in\square$
   by the induction hypothesis.

  To finish the induction step, it remains to prove $A_k\in\square$.
  As the characteristic of $\Q$ is zero, $\Q(\sqrt{A_s})$ is a Galois extension of $\Q$
  for any $s=1,\ldots,k$. Thus
  $$\Q(\sqrt{A_1},\ldots,\sqrt{A_k})=\Q(\sqrt{A_1})\cdots\Q(\sqrt{A_k})$$
  is also a Galois extension of $\Q$ in view of \cite[p.\,50, Problem 10(d)]{Mor}.
  Suppose that $A_k\not\in\square$. Then $\sqrt{A_k}\not\in\Q$,
  and hence there is an automorphism $\sigma\in\mathrm{Gal}(K/\Q)$ with $\sigma(\sqrt{A_k})\not=\sqrt{A_k}$,
  where $K=\Q(\sqrt{A_1},\ldots,\sqrt{A_k})$.
  Recall that
  $$0=x+\sum_{s=1}^k\ve_s\sqrt{A_s}W_k^{s-1}.$$
  Hence
  \begin{equation}
  \label{sigma}0=0-\sigma(0)=\sum_{s=1}^k\ve_s(\sqrt{A_s}-\sigma(\sqrt{A_s}))W_k^{s-1}.
  \end{equation}
  Note that $\sigma(\sqrt A_k)=-\sqrt{A_k}$, and $\sigma(\sqrt{A_s})\in\{\pm\sqrt{A_s}\}$
  for all $s=1,\ldots,k-1$. Thus, by \eqref{sigma} we have
  \begin{align*}2|\sqrt{A_k}|W_k^{k-1}=|2\ve_k\sqrt{A_k}W_k^{k-1}|
  \ls\sum_{s=1}^{k-1}2|\sqrt{A_s}|W_k^{s-1}.
  \end{align*}
  On the other hand,
  \begin{align*}|\sqrt{A_k}|W_k^{k-1}\gs&W_k^{k-2}\l(1+\sum_{s=1}^k|\sqrt{A_s}|\r)
  \\>&W_k^{k-2}\sum_{s=1}^{k-1}|\sqrt{A_s}|\gs\sum_{s=1}^{k-1}|\sqrt{A_s}|W_k^{s-1}.
  \end{align*}
  So we get a contradiction and this concludes our proof of the claim.

  Note that
  \begin{align*}W:=&\l(\sum_{s=1}^k(1+A_s^2)\r)\l(1+\sum_{s=1}^kA_s^{-2}\r)
  \\=&\sum_{s=1}^k(1+A_s^2)+\sum_{r=1}^k\sum_{s=1}^k A_r^{-2}(1+A_s^2).
  \end{align*}
  For $0\ls\al\ls 1$ clearly $1+\al^4\gs1\gs\al$; if $\al\gs1$ then $1+\al^4\gs\al^4\gs\al$.
  So $1+\al^4\gs \al$ for all $\al\gs0$, and hence $1+A_s^2\gs|\sqrt{A_s}|$ for all $s=1,\ldots,k$.
  Therefore,
  $$W\gs\sum_{s=1}^k(1+A_s^2)+1\gs1+\sum_{s=1}^k|\sqrt{A_s}|.$$

  If $t\in\{1,\ldots,k\}$ and $|A_t|\gs1$, then
  $$|\sqrt{A_t}|W\gs W\gs1+\sum_{s=1}^k|\sqrt{A_s}|.$$
  If $1\ls t\ls k$ and $|A_t|<1$, then  $|\sqrt{A_t}|=|A_t|^{1/2}>A_t^2$ and hence
  \begin{align*}|\sqrt{A_t}|W\gs&|\sqrt{A_t}|\l(1+\sum_{s=1}^kA_t^{-2}(1+A_s^2)\r)
  \\\gs&|\sqrt{A_t}|+\sum_{s=1}^k(1+A_s^2)=|\sqrt{A_t}|+(1+A_t^2)+\sum_{s=1\atop s\not=t}^k(1+A_s^2)
  \\\gs&1+\sum_{s=1}^k|\sqrt{A_s}|.
  \end{align*}
  Therefore the inequality \eqref{W-bound} holds if we take $W_k=W$. Applying the proved claim we immediately obtain the desired result. This concludes our proof of Theorem \ref{Th1.2}. \qed

\section{Proof of Theorem \ref{Th1.1}}
 \setcounter{equation}{0}
 \setcounter{conjecture}{0}
 \setcounter{theorem}{0}

  Let $p$ be any prime. As usual, we let $\Q_p$ and $\Z_p$ denote the $p$-adic field and the ring of $p$-adic integers respectively. We also define
  $$\Z_{(p)}=\Q\cap \Z_p=\l\{\f ab:\ a,b\in\Z\ \t{and}\ p\nmid b\r\}.$$

D. Flath and S. Wagon \cite{FW} attributed the following lemma as an observation of J. Robinson, but we cannot find it in any of Robinson's papers.

 \begin{lemma} \label{Lem3.1} Let $r$ be any rational number. Then
 \begin{equation}
 r\in\Z_{(2)}\iff \exists x\exists y\exists z[7r^2+2=x^2+y^2+z^2].
 \end{equation}
 \end{lemma}
 \Proof. The Gauss-Legendre theorem on sums of three squares (cf. \cite[pp.\, 17-23]{N96})) states that
 $n\in\N=\{0,1,\ldots\}$ is a sum of three integer squares if and only if $n\not\in\{4^k(8m+7):\ k,m\in\N\}.$

 If $r=a/b$ with $a,b\in\Z$ and $2\nmid b$, then $7a^2+2b^2\eq 2-a^2\eq1,2\pmod4$
 and hence $7a^2+2b^2$ is a sum of three squares, thus $7r^2+2=(7a^2+2b^2)/b^2$ can be expressed as $x^2+y^2+z^2$ with $x,y,z\in\Q$.

 Suppose that $r=a/b$ with $a,b\in\Z$, $2\nmid a$, $b\not=0$ and $2\mid b$. If $7r^2+2=x^2+y^2+z^2$
 for some $x,y,z\in\Q$, then there is a nonzero integer $c$ such that $c^2(7r^2+2)$ is a sum of three integer squares and hence $c^2(7r^2+2)\not\in\{4^k(8m+7):\ k,m\in\N\}$.
 Note that any odd square is congruent to $1$ modulo $8$
 and  $7a^2+2b^2\eq 7\pmod 8$ as $2\nmid a$ and $2\mid b$.
 Thus the integer
 $c^2(7r^2+2)= (c/b)^2(7a^2+2b^2)$
 has the form $(2^k)^2(8m+7)$
 with $k,m\in\N$ which leads to a contradiction.

 In view of the above, we have completed the proof of Lemma \ref{Lem3.1}.

 For any prime $p$ and $t\in\Q$, as usual we denote the $p$-adic valuation of $t$
 by $\nu_p(t)$. For $A\se\Q$ we define $A^\times=\{a\in A\sm\{0\}:\ a^{-1}\in A\}.$

\begin{lemma}\label{Lem3.2}  Let $p$ be a prime, and let $t\in\Q$. Then
\begin{equation}t\in\Z_{(p)}^\times\iff t\not=0\land (t+t^{-1}\in\Z_{(p)}).
\end{equation}
\end{lemma}
\Proof. For $t\in\Q^*$, we have $\nu_p(t^{-1})=-\nu_p(t)$. So the desired result follows. \qed

\begin{remark} This easy lemma was used by Daans \cite{D18}.
\end{remark}

For first-order formulas $\psi_{1},\ldots,\psi_k$, we simply write
$$\psi_1\lor\cdots\lor \psi_k\ \ \t{and}\ \ \psi_1\land\cdots\land\psi_k$$
as $\bigvee_{s=1}^k\psi_s$ and $\bigwedge_{s=1}^k\psi_s$ respectively.

\begin{definition} We set $\square^*=\{x^2:\ x\in\Q^*\}$.
A subset $T$ of $\Q$ is said to be {\it $m$-good} if there are polynomials
$$f_s(t,x_1,\ldots,x_m),\ g_{s1}(t,x_1,\ldots,x_m),\ldots,g_{s\ell_s}(t,x_1,\ldots,x_m)\ (s=1,\ldots,k)$$
with integer coefficients such that a rational number $t$ belongs to $T$ if and only if
$$\exists x_1\cdots\exists x_m\bigg[\bigvee_{s=1}^k\(f_s(t,x_1,\ldots,x_m)=0\land\bigwedge_{j=1}^{\ell_s}(g_{sj}(t,x_1,\ldots,x_m)
\in\square^*)\)\bigg].$$
\end{definition}
\begin{remark} \label{Rem3.1} (i) Clearly a rational number $t$ is nonzero if and only if $t^2\in\square^*$. For any polynomial $P(x)\in\Z[x]$ of degree $d$, we have $x^{2d}P(x^{-1})\in\Z[x]$, and
$$t^{2d}P(t^{-1})\in\square^*\iff P(t^{-1})\in\square^*$$
for all $t\in\Q^*$.

(ii) For any $a,b\in\Q$, clearly $(a=0\land b=0)\iff a^2+b^2=0$.
 In view of this and the distributive law concerning disjunction and conjunction, if $S\se \Q$ is $m$-good and $T\se \Q$ is $n$-good then $S\cap T$ is $(m+n)$-good.
\end{remark}

\begin{lemma}\label{LemZ2}
Both $\Z_{(2)}$ and $\Z_{(2)}^\times$ are $2$-good.
\end{lemma}
\Proof. For any $t\in\Q$, by Lemma \ref{Lem3.1} we have
$$t\in\Z_{(2)}\iff \exists x\exists y[7t^2+2-x^2-y^2\in\square].$$
Note also that
$$t\in\Z_{(2)}^\times\iff t\not=0\land (t+t^{-1}\in\Z_{(2)})$$
by Lemma \ref{Lem3.2}. Combining these with Remark \ref{Rem3.1} we immediately get the desired result. \qed

Let $a,b\in\Q^*$. As in Poonen \cite{Po}, we define
\begin{equation}S_{a,b}=\{2x_1\in\Q:\ \exists x_2\exists x_3\exists x_4[x_1^2-ax_2^2-bx_3^2+abx_4^2=1]\}
\end{equation}
and
\begin{equation}T_{a,b}=\{x+y:\ x,y\in S_{a,b}\}.\end{equation}

\begin{lemma}\label{Lem3.3} Let $a,b\in\Q^*$ with $a>0$ or $b>0$. Then
$T_{a,b}$ and $T_{a,b}^\times$ are $5$-good.
\end{lemma}
\Proof. Let $r\in\Q$. Note that
$$\l(\f r2\r)^2-a\l(\f x2\r)^2-b\l(\f y2\r)^2+ab\l(\f z2\r)^2=1\iff ab(4-r^2+ax^2+by^2)=(abz)^2.$$
So
\begin{align*}r\in S_{a,b}\iff &\exists x\exists y[ab(4-r^2+ax^2+by^2)\in\square]
\\\iff&\exists x\exists y[ab(4-r^2+ax^2+by^2)=0\lor ab(4-r^2+ax^2+by^2)\in\square^*]
\end{align*}
and hence $S_{a,b}$ is $2$-good.

For $t\in\Q$, we obviously have
$$t\in T_{a,b}\iff\exists r(r\in S_{a,b}\land t-r\in S_{a,b}).$$
As $S_{a,b}$ is $2$-good, $T_{a,b}$ is $5$-good by Remark \ref{Rem3.1}(ii).

By Koenigsmann \cite[Proposition 6]{K},
$$T_{a,b}^\times=\bigcap_{p\in\Delta_{a,b}}\Z_{(p)}^\times,$$
where $$\Delta_{a,b}=\{p:\ p\ \t{is prime and}\ (a,b)_p=-1\}$$
with $(a,b)_p$ the Hilbert symbol. (We view an empty intersection of subsets of $\Q$ as $\Q$, thus $T_{a,b}^\times=\Q$ if $\Delta_{a,b}=\emptyset$.)
Let $t\in\Q^*$. By Lemma \ref{Lem3.2}, we have
\begin{align*}t\in T_{a,b}^\times\iff \forall p\in \Delta_{a,b}(t+t^{-1}\in\Z_{(p)})
\iff t+t^{-1}\in T_{a,b}.
\end{align*}
In view of Remark \ref{Rem3.1}, from the above we see that $T_{a,b}^\times$ is $5$-good.

The proof of Lemma \ref{Lem3.3} is now complete. \qed

 For $S,T\se\Q$ we set
$$ST=\{st:\ s\in S\ \t{and}\ t\in T\}.$$
For $a,b,c\in\Q^*$ with $a>0$ or $b>0$, we define
\begin{equation}
J_{a,b}^c=T_{a,b}\,\{cy^2:\ y\in\Q\ \t{and}\ 1-cy^2\in\square T_{a,b}^\times\}.
\end{equation}
By Koenigsmann \cite[Proposition 6]{K} and Daans \cite[Lemma 5.4]{D18},
\begin{equation}\label {Jabc}
J_{a,b}^c=\bigcap_{p\in\Delta_{a,b}\atop 2\nmid \nu_p(c)}p\Z_{(p)}.
\end{equation}

\begin{lemma}\label{Lem3.4} Let $a,b,c\in\Q^*$ with $a>0$ or $b>0$. Then $J_{a,b}^c$ is $12$-good.
\end{lemma}
\Proof. As $0\in J_{a,b}^c$ by \eqref{Jabc}, we have $T_{a,b}^\times\not=\emptyset$.
For any $x\in\Q$, clearly
$$x\in\square T_{a,b}^\times\iff x=0\lor \exists y (xy^2\in T_{a,b}^\times).$$
So $\square T_{a,b}^\times$ is $6$-good in light of Lemma \ref{Lem3.3}.
As $\pm2\in S_{a,b}$, both $T_{a,b}$ and $J_{a,b}^c$ contain $0$. Let $x\in\Q$. Note that
$$x\in J_{a,b}^c\iff x=0\lor\exists y\not=0\left[\f x{cy^2}\in T_{a,b}\land (1-cy^2\in\square T_{a,b}^\times)\right].$$
Thus, with the aid of Remark \ref{Rem3.1} and Lemma \ref{Lem3.3}, we see that $J_{a,b}^c$
is $12$-good. \qed

\medskip
\noindent{\it Proof of Theorem \ref{Th1.1}}. Let $t\in\Q$. Clearly,
$$t\in \Q\sm\Z\iff t\not=0\land t^{-1}\in\bigcup_{p\in \mathbb P}p\Z_{(p)},$$
where $\mathbb P$ is the set of all primes. By Daans \cite[(1)]{D18}, we have
\begin{equation}\bigcup_{p\in\mathbb P}p\Z_{(p)}= 2\Z_{(2)}\cup\bigcup_{(a,b)\in \Phi}(J_{a,b}^a\cap J_{a,b}^{2b}),
\end{equation} where
\begin{equation}\Phi=\{(1+4u^2,2v):\ u,v\in\Z_{(2)}^\times\}.
\end{equation}
In view of this and Lemma \ref{Lem3.1}, when $t\not=0$ we have
\begin{align*}t\not\in\Z\iff& \f1{2t}\in\Z_{(2)}\lor \exists u\exists v\l[u,v\in\Z_{(2)}^\times\land \f1t\in J_{1+4u^2,2v}^{1+4u^2}\cap  J_{1+4u^2,2v}^{4v}\r]
\\\iff&\exists u\exists v\l(\f 7{4t^2}+2-u^2-v^2\in\square\r)
\\&\lor \exists u\exists v\l[u,v\in\Z_{(2)}^\times\land \f1t\in J_{1+4u^2,2v}^{1+4u^2}\cap  J_{1+4u^2,2v}^{4v}\r]
\\\iff&\exists u\exists v\bigg[8t^2+7-u^2-v^2\in\square
\\&\ \ \lor \l(u,v\in\Z_{(2)}^\times\land t^{-1}\in J_{1+4u^2,2v}^{1+4u^2}
\land t^{-1} \in J_{1+4u^2,2v}^{4v}\r)\bigg].
\end{align*}
Combining this with Lemmas \ref{LemZ2} and \ref{Lem3.4},
 we obtain that $\Q\sm\Z$ is $30$-good in view of Remark \ref{Rem3.1}.

By the above, there are polynomials
$$f_s(t,x_1,\ldots,x_{30}),\ g_{s1}(t,x_1,\ldots,x_{30}),\ldots,g_{s\ell_s}(t,x_1,\ldots,x_{30})
\ (s=1,\ldots,k)$$
with integer coefficients such that a rational number $t$ is not an integer if and only if
$$\exists x_1\cdots\exists x_{30}\bigg[\bigvee_{s=1}^k\(f_s(t,x_1,\ldots,x_{30})=0\land\bigwedge_{j=1}^{\ell_s}(g_{sj}(t,x_1,\ldots,x_{30})
\in\square^*)\)\bigg].$$
Note that
\begin{equation*}\label{g}g_{sj}(t,x_1,\ldots,x_{30})\not=0\quad\t{for all}\ j=1,\ldots,\ell_s
\end{equation*}
if and only if
$$x_{31}\prod_{j=1}^{\ell_s}g_{sj}(t,x_1,\ldots,x_{30})-1=0$$
for some $x_{31}\in\Q$. By Theorem \ref{Th1.2}, when $\prod_{j=1}^{\ell_s}g_{sj}(t,x_1,\ldots,x_{30})\not=0$, we have
$$g_{sj}(t,x_1,\ldots,x_{30})\in\square\qquad\t{for all}\ j=1,\ldots,\ell_s$$
if and only if
$${\mathcal J}_{\ell_s}\l(g_{s1}(t,x_1,\ldots,x_{30}),\ldots,
g_{s\ell_s}(t,x_1,\ldots,x_{30}),x_{32}\r)=0$$
for some $x_{32}\in\Q$. Combining these we see that
$t\not\in\Z$
if and only if there are $x_1,\ldots,x_{32}\in\Q$ such that the product of all those
\begin{align*}&f_s(t,x_1,\ldots,x_{30})^2+\(x_{31}\prod_{j=1}^{\ell_s}g_{sj}(t,x_1,\ldots,x_{30})-1\)^2
\\&+{\mathcal J}_{\ell_s}(g_{s1}(t,x_1,\ldots,x_{30}),\ldots,g_{s\ell_s}(t,x_1,\ldots,x_{30}),x_{32})^2
\end{align*}
$(s=1,\ldots,k)$ is zero.

In the spirit of the above proof, we can actually construct an explicit polynomial $P(t,x_1,\ldots,x_{32})$
with integer coefficients satisfying \eqref{32} with the total degree of $P$ smaller than $2.1\times10^{11}$. This concludes our proof of Theorem \ref{Th1.1}. \qed

\section{Proof of Theorem \ref{Th1.3}}
 \setcounter{equation}{0}
 \setcounter{conjecture}{0}
 \setcounter{theorem}{0}

  It is known that each nonnegative integer can be written as a sum of four squares of rational numbers.
  This result due to Euler (cf. \cite{Pi}) is weaker than Lagrange's four-square theorem (cf. \cite[pp.\,5-7]{N96}).
   Clearly, any nonnegative rational number can be written as $a/b=(ab)/b^2$ with $a,b\in\N$ and $b>0$.
  So we have the following lemma.

 \begin{lemma}\label{Lem4.1} Let $r\in\Q$. Then
 \begin{equation}r\gs0\iff\exists w\exists x\exists y\exists z[r=w^2+x^2+y^2+z^2].
 \end{equation}
 \end{lemma}

 We also need a known result of Sun \cite[Theorem 1.1]{S21}.

 \begin{lemma}[Sun \cite{S21}]\label{Lem4.2}
 Let $\mathcal A\se\N$ be an r.e. (recursively enumerable) set.

  {\rm (i)} There is a polynomial  $P_{\mathcal A}(x_0,x_1,\ldots,x_{9})$ with integer coefficients such that for any $a\in \N$ we have
 $a\in\mathcal A$ if and only if $P_{\mathcal A}(a,x_1,\ldots,x_{9})=0$
 for some $x_1,\ldots,x_{9}\in\Z$ with $x_{9}\gs0$.

 {\rm (ii)} There is a polynomial  $Q_{\mathcal A}(x_0,x_1,\ldots,x_{10})$ with integer coefficients such that for any $a\in \N$ we have
 $a\in\mathcal A$ if and only if $Q_{\mathcal A}(a,x_1,\ldots,x_{10})=0$
 for some $x_1,\ldots,x_{10}\in\Z$ with $x_{10}\not=0$.
 \end{lemma}

 \medskip
 \noindent{\it Proof of Theorem \ref{Th1.3}}. It is well known that there are nonrecursive r.e. sets
 (see, e.g., \cite[pp.\,140-141]{C}). Let us take any nonrecursive r.e. set $\mathcal A\se\N$.

 {\rm (i)} Let $P_{\mathcal A}$ and $P$ be polynomials as in Lemma \ref{Lem4.2} and Theorem \ref{Th1.1}. In view of  Lemmas \ref{Lem4.1}-\ref{Lem4.2} and Theorem \ref{Th1.1}, for any $a\in\N$ we have
 \begin{align*} a\not\in\mathcal A\iff& \forall x_1\cdots\forall x_{9}[\neg(x_1,\ldots,x_{9}\in\Z\land x_{9}\ge0)
 \lor P_{\mathcal A}(a,x_1,\ldots,x_{9})\not=0]
 \\\iff&\forall x_1\cdots\forall x_{9}\bigg[\bigvee _{t=1}^{9}(x_t\not\in\Z)\lor x_{9}<0
 \lor P_{\mathcal A}(a,x_1,\ldots,x_{9})\not=0\bigg]
 \\\iff&\forall x_1\cdots\forall x_{9}\bigg[\bigvee _{t=1}^{9}\exists y_1\cdots\exists y_{32}(P(x_t,y_1,\ldots,y_{32})=0)
 \\&\lor -x_{9}>0\lor \exists y_1(y_1P_{\mathcal A}(a,x_1,\ldots,x_{9})-1=0)\bigg]
 \\\iff&\forall x_1\cdots\forall x_{9}\exists y_1\cdots\exists y_{32}
 [P_0(a,x_1,\ldots,x_9,y_1,\ldots,y_{32})=0],
 \end{align*}
 where
 \begin{align*}&P_0(a,x_1,\ldots,x_9,y_1,\ldots,y_{32})
 \\=&(y_1P_{\mathcal A}(a,x_1,\ldots,x_{9})-1)\prod_{t=1}^{9}P(x_t,y_1,\ldots,y_{32})
 \\&\times\l((x_{9}y_1-1)^2+(x_9+y_2^2+y_3^2+y_4^2+y_5^2)^2\r).
 \end{align*}
 It follows that for any $a\in\N$ we have
 $$a\in\mathcal A\iff \exists x_1\cdots\exists x_9\forall y_1\cdots\forall y_{32}
 \exists y_{33}[y_{33}P_0(a,x_1,\ldots,x_9,y_1,\ldots,y_{32})-1=0]
 $$
 As both $\mathcal A$ and $\N\sm\mathcal A$ are nonrecursive, by the above we get that $\forall^9\exists^{32}$ over $\Q$ and $\exists^9\forall^{32}\exists$ over $\Q$ are undecidable.

  {\rm (ii)} Let $Q_{\mathcal A}$ be the polynomial in Lemma \ref{Lem4.2}(ii). For any $a\in\N$, we have
  \begin{align*} a\not\in\mathcal A\iff& \forall x_1\cdots\forall x_{10}[\neg(x_1,\ldots,x_{10}\in\Z\land x_{10}\not=0)
 \lor Q_{\mathcal A}(a,x_1,\ldots,x_{10})\not=0]
 \\\iff&\forall x_1\cdots\forall x_{10}\bigg[\bigvee_{t=1}^{10}(x_t\not\in\Z)\lor x_{10}=0
 \lor Q_{\mathcal A}(a,x_1,\ldots,x_{10})\not=0\bigg].
 \end{align*}
  By the proof of Theorem \ref{Th1.1}, $\Q\sm\Z$ is $30$-good.
 Thus, in view of Theorem \ref{Th1.2},
 there are polynomials
 $$f_s(x,y_1,\ldots,y_{31})\ \t{and}\ g_s(x,y_1,\ldots,y_{31})\ \ (s=1,\ldots,k)$$
 with integer coefficients such that for any $x\in\Q$ we have
 $$x\not\in\Z\iff \exists y_1\cdots\exists y_{31}\bigg[\bigvee_{s=1}^k(f_s(x,y_1,\ldots,y_{31})=0\land
 g_s(x,y_1,\ldots,y_{31})\not=0)\bigg]$$
 Thus, for any $a\in\N$, we have
 \begin{align*} a\not\in\mathcal A\iff& \forall x_1\cdots\forall x_{10}\exists y_1\cdots\exists y_{31}
 \\&\bigg[\bigvee_{t=1}^{10}\bigg(\bigvee_{s=1}^k(f_s(x_t,y_1,\ldots,y_{31})=0\land
 g_s(x_t,y_1,\ldots,y_{31})\not=0)
 \\&\lor x_{10}=0
 \lor Q_{\mathcal A}(a,x_1,\ldots,x_{10})\not=0\bigg)\bigg]
 \end{align*}
 and hence
 \begin{align*} a\in\mathcal A\iff& \exists x_1\cdots\exists x_{10}\forall y_1\cdots\forall y_{31}
 \\&\bigg[\bigwedge_{t=1}^{10}\bigg(\bigwedge_{s=1}^k(f_s(x_t,y_1,\ldots,y_{31})\not=0\lor
 g_s(x_t,y_1,\ldots,y_{31})=0)
 \\&\land x_{10}\not=0
 \land Q_{\mathcal A}(a,x_1,\ldots,x_{10})=0\bigg)\bigg].
  \end{align*}
  Let $\Gamma=\{1,\ldots,k\}\times\{1,\ldots,10\}$. By the distributive law concerning disjunction and conjunction,
  $$\bigwedge_{t=1}^{10}\bigwedge_{s=1}^k\l(f_s(x_t,y_1,\ldots,y_{31})\not=0\lor
 g_s(x_t,y_1,\ldots,y_{31})=0\r)$$
 is equivalent to
 \begin{align*}\bigvee_{\Delta\se \Gamma}
 \bigg(\bigwedge_{(s,t)\in\Delta}(f_s(x_t,y_1,\ldots,y_{31})\not=0)
 \land\bigwedge_{(s',t')\in\Gamma\sm\Delta}(g_{s'}(x_{t'},y_1,\ldots,y_{31})=0)\bigg).
 \end{align*}
 Thus, for any $a\in\N$, we have
 \begin{align*}a\in \mathcal A
 \iff&\exists x_1\cdots\exists x_{10}\forall y_1\cdots\forall y_{31}
 \\&\bigg[\bigvee_{\Delta\se\Gamma}\bigg(x_{10}\prod_{(s,t)\in\Delta}f_s(x_t,y_1,\ldots,y_{31})\not=0
 \\&\land\bigwedge_{(s',t')\in\Gamma\sm\Delta}(g_{s'}(x_{t'},y_1,\ldots,y_{31})=0)
  \land Q_{\mathcal A}(a,x_1,\ldots,x_{10})=0\bigg)\bigg]
 \\\iff&\exists x_1\cdots\exists x_{10}\forall y_1\cdots\forall y_{31}\exists z
 \\&\bigg[\bigvee_{\Delta\se \Gamma}\bigg(1-zx_{10}\prod_{(s,t)\in\Delta}f_s(x_t,y_1,\ldots,y_{31})=0
 \\&\land \bigwedge_{(s',t')\in\Gamma\sm\Delta}(g_{s'}(x_{t'},y_1,\ldots,y_{31})=0)
 \land Q_{\mathcal A}(a,x_1,\ldots,x_{10})=0\bigg)\bigg]
 \end{align*}
 and hence
 \begin{align*}a\in \mathcal A
 \iff&\exists x_1\cdots\exists x_{10}\forall y_1\cdots\forall y_{31}\exists z
 [P_1(a,x_1,\ldots,x_{10},y_1,\ldots,y_{31},z)=0],
 \end{align*}
 where we view an empty product as $1$, and $P_1(a,x_1,\ldots,x_{10},y_1,\ldots,y_{31},z)$ stands for the product of
\begin{align*}&\bigg(1-zx_{10}\prod_{(s,t)\in\Delta}f_s(x_t,y_1,\ldots,y_{31})\bigg)^2
 \\+&\sum_{(s',t')\in\Gamma\sm\Delta}g_{s'}(x_{t'},y_1,\ldots,y_{31})^2+Q_{\mathcal A}(a,x_1,\ldots,x_{10})^2
 \end{align*}
 over $\Delta\se\Gamma$.
As $\mathcal A$ is nonrecursive, we obtain that $\exists^{10}\forall^{31}\exists$ over $\Q$
is undecidable.

In view of the above, we have completed the proof of Theorem \ref{Th1.3}. \qed

\Ack. The authors would like to thank the referee for helpful comments.

 \end{document}